\documentclass[runningheads]{llncs}

\usepackage{graphicx}
\usepackage{amsmath}
\usepackage{amssymb}
\usepackage{hyperref}
\usepackage[capitalize]{cleveref}
\usepackage{subfigure}
\usepackage{verbatim}

\crefformat{equation}{(#2#1#3)}
\makeatletter
\newcommand*{\textoverline}[1]{$\overline{\hbox{#1}}\m@th$}
\makeatother

\begin{document}
\title{Algorithm~and~abstraction~in~formal~mathematics}

\author{Heather Macbeth\orcidID{0000-0002-0290-4172}}

\institute{Fordham University, New York NY 10023, USA\\
\email{hmacbeth1@fordham.edu}}
\maketitle
\begin{abstract}
  I analyse differences in style between traditional prose mathematics writing and computer-formalised mathematics writing,
  presenting five case studies.
  I note two aspects where good style seems to differ between the two: in
  their incorporation of computation and of abstraction.
  I argue that this reflects a different mathematical aesthetic for formalised mathematics.
\end{abstract}
\section{Introduction}

In the last twenty years, formalisation---%
building up proofs as line-by-line logical deductions from the axioms of mathematics, with the help of specialised computer systems\footnote{
  Examples include Agda, Coq, Lean, HOL Light, Isabelle, Metamath and Mizar.}---%
has seen increasing interest from mathematicians. 
The rapidly increasing coverage of the mathematical literature in these systems
is very much a social process:
their mathematical libraries are built collaboratively by hundreds of people,
and code contributed by one person will be reviewed in detail by another,
and often thoroughly re-worked a year later by a third.

In this kind of human and social process, culture develops spontaneously.
The back-and-forth of discussion in this process includes frequent comment
on a formalised proof's beauty, elegance, cleverness,
and other abstract properties generally associated with mathematical aesthetics.
The communities of mathematicians doing this work
consider computer-formalised proofs to be,
not simply utilitarian certificates for the correctness of logical claims,
but a fully-fledged medium for mathematical exposition.

In this article I describe
(necessarily very subjectively)
some aspects of this aesthetic of computer-formalised mathematics writing.
Much of this aesthetic is inherited from traditional prose mathematics writing,
on which there is a vast literature \cite{Arn98,Bon82,Hal70,Har40,IA14,Mon12,Rot97,Tao07,Wel90}.
I therefore focus on cases in which good style in formalised mathematics
seems to differ from good style in traditional prose mathematics.
I present five case studies,\footnote{
Disproportionately drawn from Lean's \cite{Lean3,Lean4} 
\emph{Mathlib} \cite{Mathlib},
of which I am a maintainer.}
grouped by theme:
how to integrate computation (\cref{sec:computation})
and how much use to make of abstraction (\cref{sec:abstraction}).

\section{Computation} \label{sec:computation}

A faithful computer-formalised translation of a traditional prose proof will commonly use computation
``in the small:''
a proof step which seems obvious to humans
often represents a whole chain of strict logical reasoning,
and in most modern systems automation is used to help construct such chains.

Interestingly, such a translation will sometimes also use computation ``in the large:''
several notable formalisations' \cite{Gon08Notices,Hal17,Imm2018}
targets are
theorems whose published proofs rely on the result reported by a computer program.

So what about using computation ``in the middle?"
In this section I explore proofs where there is no absolute need to outsource a calculation to computer---%
and where, in traditional writing, simple inertia would prevent one from doing this---%
but which are arguably improved by increased reliance on computation.

\subsection{Classification of wallpaper groups} \label{sec:wallpaper}

My first example 
arises in classifying the 17 wallpaper groups.
This classification is heavily dependent on case analysis,
one branch of which is to consider wallpaper groups which contain translations and
rotations but no reflections.
These can be classified according to the orbit types of centres of symmetry.
For example, one of these wallpaper groups,
which in our classification we will associate to the tuple $(2,3,6)$,
has three centres of symmetry,
at which the stabilisers are generated by rotations of
$\frac{2\pi}{2}$,
$\frac{2\pi}{3}$,
and $\frac{2\pi}{6}$.

\begin{figure}[!htb]
  \centering
  \subfigure[$(3,3,3)$]{\includegraphics[width=0.3\textwidth, trim = 0 0 180 180, clip]{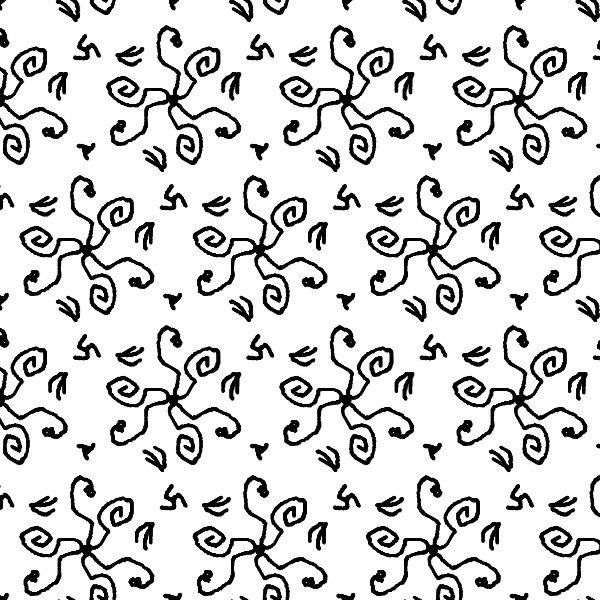}}\hfill
  \subfigure[$(2,4,4)$]{\includegraphics[width=0.3\textwidth, trim = 0 0 180 180, clip]{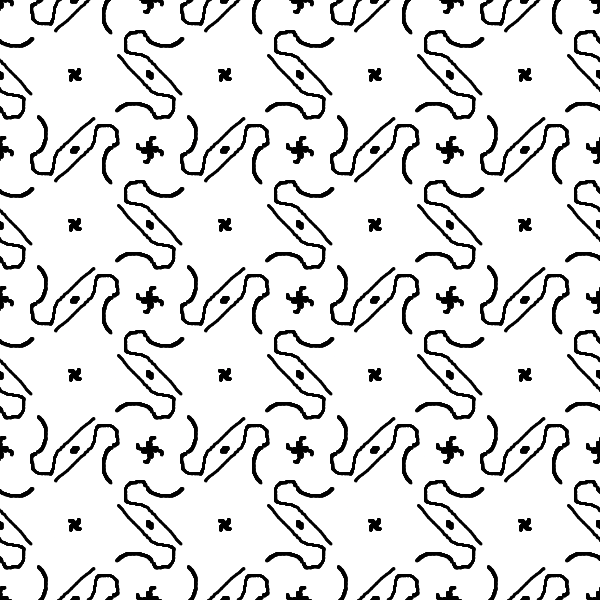}}\hfill
  \subfigure[$(2,3,6)$]{\includegraphics[width=0.3\textwidth, trim = 0 0 180 180, clip]{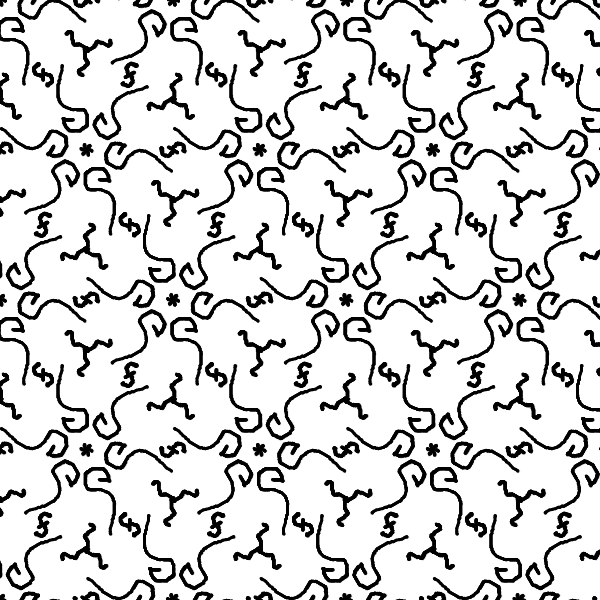}}
  \caption{The wallpaper groups containing
  rotations but no reflections. (Images via \cite{Eck}.)}
  \label{fig:wallpaper}
\end{figure} 

The following arithmetic lemma
classifies the possible tuples which can arise.
The wallpaper groups associated to these tuples are depicted in \cref{fig:wallpaper}.

\begin{lemma}
  Let $2\le p \le q \le r$ be natural numbers, with
  \begin{equation}
    \frac{1}{p}+\frac{1}{q}+\frac{1}{r}=1. \tag{$\star$}
  \end{equation}
  Then $(p, q, r)$ is one of $(3,3,3)$, $(2,4,4)$, $(2,3,6)$.
\end{lemma}

I first present a traditional prose proof
lifted from a textbook \cite{CBG08}.
\begin{proof}
  We get $p=q=r=3$ if all of $\frac{1}{p}$, $\frac{1}{q}$ and $\frac{1}{r}$ have their mean value
  of $\frac{1}{3}$.  Otherwise $p$ must be 2.

  If $r$ and $q$ have \emph{their} mean value of $\frac{1}{4}$, we get $p=2$, $q=r=4$.

  If not, $q$ must be 3, and $r$ is forced to be 6, by $(\star)$.
\end{proof}

Secondly, I describe a proof of this lemma that I wrote in Lean with Anne Baanen.
I am trying to translate the Lean code fairly literally.

\begin{proof}
The inequalities
$0<\frac{1}{r}\le \frac{1}{q} \le \frac{1}{p} \le \frac{1}{2}$ and the equality $(\star)$ yield that
\begin{alignat}{3}
  \tfrac{1}{3}&\le& \tfrac{1}{p}&\le \tfrac{1}{2}\label{eq:p}\\
  \tfrac{1}{2}\left[1-\tfrac{1}{p}\right]&\le &\tfrac{1}{q}&<\tfrac{1}{2}\label{eq:q}\\
  &&\tfrac{1}{r}&=1-\tfrac{1}{p}-\tfrac{1}{q}.\label{eq:r}
\end{alignat}
There are finitely many natural numbers $p$ satisfying \cref{eq:p}; case-split on these.
For each of these there are finitely many natural numbers $q$ satisfying \cref{eq:q}; case-split on these.
For each of these, $r$ can be determined from \cref{eq:r}.
\end{proof}

There is an algorithm implicit in these proofs.
The second (formalised) proof looks almost like a recipe for cooking the first (textbook) proof:
it describes the steps to be carried out,
rather than actually performing those steps visibly for the reader
(i.e.\ documenting the available choices at each case split).

This is very typical:
as mentioned, proof-writing in systems such as Lean frequently invokes ``tactics,''
small computer programs to construct parts of proofs.
But once we start to describe proofs via the recipes which would construct them,
there is no need to stick to the original recipe.
This was noted by Hales et~al.~\cite{Hal17}:
\begin{quote}
In the original, computer calculations were a last resort after as much was done by
hand as feasible. In the [formalisation], the use of computer has been fully embraced.
As a result, many laborious lemmas of the original proof can be automated or
eliminated altogether.
\end{quote}
I will argue that an aesthetically pleasing formal proof is one which has a short and simple recipe.
As the next two examples will show,
this is not the same thing as a proof which is itself short and simple.

\subsection{The Kochen--Specker paradox}

I next consider
a theorem from quantum mechanics.
\begin{theorem}[Kochen--Specker \cite{KS67}]
  There does not exist a boolean (say red or green) colouring of the vectors in $\mathbb{R}^3$,
  such that all triples $u,v,w\in \mathbb{R}^3$ of nonzero mutually-orthogonal vectors
  are coloured green, red, red in some order.
\end{theorem}

I will discuss a streamlined proof due to Peres \cite{Per91}.
The approach is to deduce a contradiction from the colouring of the following 33
  nonzero vectors\footnote{
      Down from 117 vectors in the original Kochen--Specker proof.
        Following Peres' notation, \textoverline{1} is shorthand for $-1$, 2 is shorthand for $\sqrt{2}$, and \textoverline{2} is
shorthand for $-\sqrt{2}$.
  } in
$\mathbb{R}^3$:

\begin{align}
  \begin{tabular}{|c|c|c|c|c|c|c|c|c|c|c|}
  \textoverline{1}\textoverline{1}2 & \textoverline{1}02 &
  \textoverline{1}12 & \textoverline{1}2\textoverline{1} &
  \textoverline{1}20 & \textoverline{1}21 &
  0\textoverline{1}2 & 002 &
  012 & 02\textoverline{2} &
  02\textoverline{1} \\ 020 &
  021 & 022 &
  1\textoverline{1}2 & 102 &
  112 & 12\textoverline{1} &
  120 & 121 &
  2\textoverline{2}0 & 2\textoverline{1}\textoverline{1} \\
  2\textoverline{1}0 & 2\textoverline{1}1 &
  20\textoverline{2} & 20\textoverline{1} &
  200 & 201 &
  202 & 21\textoverline{1} &
  210 & 211 &
  220
\end{tabular}\label{KS-witnesses}
\end{align}
The basic driver of the proof is that,
once enough of the 33 vectors have been coloured,
the colours of the rest can be determined greedily:
a vector orthogonal to a green vector must be red,
and a vector orthogonal to two orthogonal red vectors must be green.

\begin{figure}
  \centering
\includegraphics*[scale=0.25]{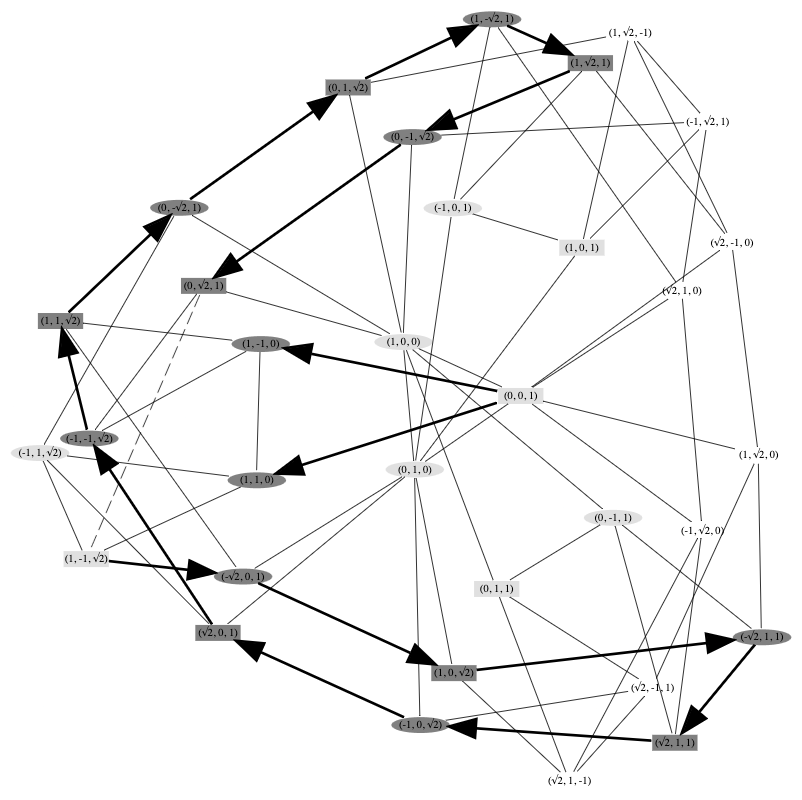}
\caption{Peres' argument \cite{Per91} for the impossibility of a suitable boolean colouring of the 33 vectors \cref{KS-witnesses}.
Edges connect orthogonal vectors.
Oval (respectively, rectangular) labels denote vectors found to be red (resp.\ green).
Light-shaded (resp.\ dark-shaded) labels denote vectors whose colour is chosen by an initial symmetry argument
(resp.\ is forced as a consequence of adjacent vectors' colours).
Arrows indicate the order in which these deductions occur.
The dotted edge indicates two orthogonal vectors which are both green, contradiction.}
\label{fig:peres-graph}
\end{figure}

Here is an outline of Peres' proof \cite{Per91} of the impossibility.
\begin{proof}
  We can determine the colours of some vectors without loss of generality.
  
{\footnotesize
  \begin{tabular}{l|l|l}
  By the symmetry
  & and by the known facts
  & we can assume\\
  \hline
  choice of $z$-axis 
  && $001$ green;  $100$, $010$ red\\
  choice of $x$ vs $-x$ 
  &$010$ red
  & $101$ green; $\overline{1}01$ red\\
  choice of $y$ vs $-y$ 
  &$100$ red
  & $011$ green; $0\overline{1}1$ red\\
  choice of $x$ vs $y$ 
  &$001$ green, thus $110$ red
  & $1\overline{1}2$ green; $\overline{1}12$ red
\end{tabular}}

Now a suitable greedy sequence of deductions 
[depicted in \cref{fig:peres-graph},
written out explicitly in the original]
forces a contradiction.
\end{proof}

For comparison, here is an outline of a proof formalised by John Harrison in HOL Light in 2005.\footnote{
  \url{https://github.com/jrh13/hol-light/blob/e736197/Tutorial/Custom_tactics.ml}}
It is really a brute force search.

\begin{proof}
Perform the following binary search:
split on a vector whose colour is not yet known;
then in each case (red or green) greedily make all possible deductions.
Stop if a contradiction is found. Recurse if not.

The result of this process is that every branch terminates in a contradiction.
\end{proof}

The most notable difference from Peres' prose proof
is to abandon his
symmetry argument (which reduces to only one configuration
on which the greedy algorithm need be run)
and instead just run the greedy argument at each stage of a binary search.\footnote{
A second difference is that when running the greedy algorithm from a partial colouring produces a contradiction,
the formalised version does not write out the \emph{certificate}:
an explicit path of deductions leading to the concluding contradiction.
But this is less controversial.
In the case of the configuration in \cref{fig:peres-graph},
when the path of deductions is written out,
it does not appear to contain any particular insight.
See Harrison \cite[section 3.4]{Har96} for a similar example.}
In effect, both arguments amount to the implementation of a search algorithm.
The search algorithm used in the formalised proof is very simple,
whereas the symmetry considerations incorporated in the original proof
can be considered as baroque optimisations to the search algorithm
to get its ``runtime" within the scale of human readability.
Harrison \cite[section 3]{Har09} defends this choice in his formalisation on the grounds of convenience:
the simpler algorithm is easier to implement.
But I argue that it is also defensible on the grounds of aesthetics:
the simpler algorithm is easier for the reader to grasp.\footnote{
Gonthier et al. \cite[section 4.3]{Gon13} report a similar example,
in which an appeal to a logical decision procedure
produces an ``intellectually more satisfying'' proof
than the original argument involving a detailed combinatorial case analysis.}

\subsection{Multiplication of Chebyshev polynomials}

My last example on the theme of computation highlights a different kind of computation.
Let $T_n(x)$ denote the $n$-th Chebyshev polynomial of the first kind.  Recall these polynomials
  satisfy a recurrence relation
  \[
    T_{n+2}(x)=2xT_{n+1}(x)-T_n(x).
  \]

\begin{lemma}[Multiplication formula for Chebyshev polynomials]
  For all natural numbers $m$ and $k$, $2T_mT_{m+k}=T_{2m+k}+T_k$.
\end{lemma}

A purely algebraic\footnote{There is an alternative approach using trigonometric identities.}
proof of this lemma is necessarily an induction.
The 
inductive step of the
scheme that works has us prove a statement for $m+2$,
\[
  \forall k :\mathbb{N}, 2T_{m+2}T_{(m+2)+k}=T_{2(m+2)+k}+T_{k},
\]
given the corresponding statements for $m$ and $m+1$.

The following is how we have been trained to write rigorous proofs of equalities in mathematics articles:
as a transitive chain of reasoning.

  \begin{proof}
Indeed,
\begin{align*}
2T_{m+2}T_{m+k+2}
&=2[2xT_{m+1}-T_m]T_{m+k+2}\\
&=2x[2T_{m+1}T_{(m+1)+(k+1)}]-2T_mT_{m+(k+2)}\\
&=2x[T_{2(m+1)+(k+1)}+T_{k+1}]-[T_{2m+(k+2)}+T_{k+2}]\\
&=[2xT_{2m+k+3}-T_{2m+k+2}]+[2xT_{k+1}-T_{k+2}]\\
&=T_{2m+k+4}+T_{k}.
\end{align*}
  \end{proof}

Halmos \cite[section 16]{Hal70},
writing long before interactive proof assistants were widespread,
 calls out the ``proof that consists of a chain of expressions separated by equal signs"
as an example of lazy writing,
\begin{quote}
unhelpful [symbolism]
that merely reports the result of the act and leaves the reader to guess how they were obtained, 
\end{quote}
and advocates for replacing such proofs by a ``recipe for action''
(a metaphor I already borrowed in \cref{sec:wallpaper}).
Here is an alternate proof of the Chebyshev lemma which precisely consists of such a recipe.
This approach follows a formalisation of mine, contributed to Mathlib.\footnote{
  Mathlib \cite{Mathlib}, \href{https://github.com/leanprover-community/mathlib4/blob/21b3a48ba0f69a4ae9eb44dbacbc4de43a6290bf/Mathlib/RingTheory/Polynomial/Chebyshev.lean\#L209}{\texttt{RingTheory/Polynomial/Chebyshev}, line 209}
}
\begin{proof}
  Indeed, two applications of the inductive hypothesis give
\begin{align}
  2T_{m+1}T_{(m+1)+(k+1)}&=T_{2(m+1)+(k+1)}+T_{k+1}\tag{$\star_1$}\\
  2T_mT_{m+(k+2)}&=T_{2m+(k+2)}+T_{k+2}\tag{$\star_2$}
\end{align}
and three applications of the recurrence relation give
\begin{align*}
T_{m+2}&=2xT_{m+1}-T_m\tag{$*_1$}\\
T_{(2m+k+2)+2}&=2xT_{(2m+k+2)+1}-T_{2m+k+2}\tag{$*_2$}\\
T_{k+2}&=2xT_{k+1}-T_{k}\tag{$*_3$}
\end{align*}
A Gr\"obner basis computation\footnote{
  This ability to send a computation to the Gr\"obner basis algorithm
  is a standard offering in formalisation languages \cite{Har07,CW07,Pot08}.
  In Lean this is performed via an external call to Sage;
  it was implemented by Dhruv Bhatia and Rob Lewis in 2022.}
shows that LHS $-$ RHS of the desired result,
\begin{align*}
2T_{m+2}T_{m+k+2}&=T_{2m+k+4}+T_{k},
\end{align*}
is in the ideal generated by LHS $-$ RHS of $(\star_1)$, $(\star_2)$, $(*_1)$, $(*_2)$, $(*_3)$.
\end{proof}

In a traditional prose proof,
there is a high barrier to outsourcing this kind of computation to a specialised computer algebra system.
The code performing the calculation must set up (under some names) the 11 variables\footnote{
  We normalise indices before the computation.
}
\[
  \begin{tabular}{llllll}
  $T_{m+2}$,&
  $T_{k+2}$,&
  $T_{m+k+2}$,&
  $T_{2m+k+4}$,&
  $x$,\\
  $T_{m+1}$,&
  $T_{k+1}$,&&
  $T_{2m+k+3}$,\\
  $T_m$,&
  $T_{k}$,&&
  $T_{2m+k+2},$
  \end{tabular}
\]
the five polynomials in these 11 variables which generate the ideal,
and a sixth polynomial whose membership in the ideal is to be checked.
This process is tedious and error-prone;
it will demand close attention from both author and reader.
By contrast, when formalising, there is no such barrier:
the problem statement is already available in a suitable electronic format.

The point is not just that in formalisation the second proof becomes feasible;
I argue  it is also more elegant.
It is easier to grasp at high level:
it is clear upfront what facts are being used,
and the reader can check by eye that the goal appears to be within the scope of the Gr\"obner basis algorithm as run on these facts.
Its black-boxing of the routine algorithm also makes the ideas more transparent---%
in this case, the choice of specialisations of the two 
inductive hypotheses.\footnote{
Indeed, let $k+b$ be the chosen instantiation of the $(m+1)$-inductive hypothesis
and $k+a$ that of the $m$-inductive hypothesis:
\begin{align*}
  2T_{m+1}T_{(m+1)+(k+b)}&=T_{2(m+1)+(k+b)}+T_{k+b},\\
  2T_mT_{m+(k+c)}&=T_{2m+(k+c)}+T_{k+c}.
\end{align*}
In order for there to be a nontrivial polynomial relation among these, the goal
\[
  2T_{m+2}T_{(m+2)+k}=T_{2(m+2)+k}+T_{k},
\]
and some uses of the recurrence,
we need to arrange that the $T$-indices which appear, 
\[
  m+\{2,1,0\},\quad k+\{0,b, c\},\quad
  m+k+\{2, b+1, c\},\quad 2m+k+\{4, b+2, c\},
\]
are all either (a) sets of three consecutive numbers
(in which case the recurrence relation provides an identity connecting them)
or (b) all the same.
This forces $c=2$, and that forces $b=1$,
leading to the instantiations $(\star_1)$, $(\star_2)$ chosen.}

In summary, I argue that in formalisation the threshold for switching to full automation 
should lower,
with many ``mid-sized" computations automated away.

\section{Abstraction} \label{sec:abstraction}

I now turn to the second realm in which I argue 
that there is a stylistic difference between prose and formal mathematics:
the question of abstraction.

The principle that every mathematical argument should be generalised to exactly its proper context
dates at least to Bourbaki \cite[section 2]{Bou50}:
\begin{quote}
Where the superficial observer sees only two, or several, quite
distinct theories, lending one another ``unexpected support"
\ldots
[we advocate]
to look for the deep-lying reasons for such a discovery, 
to find the common ideas of these theories,
buried under the accumulation of details properly belonging to each of them \ldots
and to put them in their proper light.
\end{quote}
This idea was profoundly influential.
But though widely agreed on in principle,
it is not followed universally in practice.
For example, Halmos \cite{Hal70} advises  writers,
\begin{quote}
  The observation that a proof proves something a little more general than it was invented for can frequently be left to the reader.
\end{quote}
The main reason is psychological:
abstractions seem to be a cognitive barrier for readers.
A secondary, related reason is practical:
you can't expect your reader to be confident in the application of an abstraction
that she has never seen before,
and so it's courteous to her to specialise it.

In formalised mathematics the trade-offs are different.
The practical obstruction to abstraction nearly disappears,\footnote{
Your reader has immediate access to a full exposition of an unfamiliar abstraction;
moreover, thanks to verification, she can trust you
when you state that all the preconditions hold for that abstraction to be applicable in the context at hand.}
though the psychological one remains.
Moreover, as the examples in this section will explore,
the usual arguments in favour of abstraction apply somewhat more strongly than in prose mathematics writing.
All told, formal mathematics favours decidedly more use of abstraction.

\subsection{Lax--Milgram theorem}
I first consider the Lax--Milgram theorem,
a functional analysis result
which turns up in the standard approach to linear elliptic partial differential equations.

Let $H$ be a real Hilbert space, $B : H \times H \to \mathbb{R}$ a bilinear form.

\begin{theorem}[Lax--Milgram] Suppose there exist constants $\alpha, \beta > 0$ so that
\begin{itemize}
  \item (boundedness) for all $u,v\in H$, $|B[u, v]| \le \alpha \|u\| \|v\|$
  \item (coercivity) for all $u \in H$, $B[u, u]\geq \beta\|u\|^2$.
\end{itemize}
Then for each $f \in H^*$, there exists a unique $u \in H$ so that
for all $v \in H$, $B[u, v] = f(v)$.
\end{theorem}

The proof of this theorem
begins by constructing a bounded linear map $A:H\to H$ such that for all $u,v\in H$,
we have $B[u,v]=\langle A(u),v\rangle$.
By the coercivity of $B$, we have for all $u$
\[
\beta\|u\|^2\le B[u, u]=\langle A(u), u\rangle \le \|A(u)\| \|u\|,
\]
so (by the above if $u\ne 0$ and trivially if $u = 0$)
\begin{equation}
  \beta\|u\|\le \|A(u)\|.\label{eq:A-antilipschitz}
\end{equation}

It suffices to show that the operator $A$ is bijective.
I will concentrate on one step of the bijectivity argument:
the step where we exploit \cref{eq:A-antilipschitz} to establish that $A$
is injective and has closed range.
As usual I present two proofs.

\begin{proof}[\cite{How20}, slightly compressed]
If $u_1, u_2 \in H$,
then $\lVert A(u_1 - u_2)\rVert \geq \beta \lVert u_1 - u_2\rVert$,
from which it's clear that $A$ is injective.

To see that the range of $A$ is closed in $H$,
let $\{u_j\}_{j=1}^\infty\subset H$ satisfy $Au_j \to w$ for some $w \in H$.
We need to show that there exists $u \in H$
so that $Au = w$.

For this, we notice that
\[
  \lVert u_i - u_j \rVert \le \beta^{-1} \lVert Au_i - Au_j \rVert.
\]
The sequence $\{Au_j\}_{j=1}^\infty$ converges, so it must be Cauchy,
so we see that $\{u_j\}_{j=1}^\infty$ must be Cauchy,
and so must converge to some $u \in H$.
Since $A$ is bounded,
\[
\lVert Au - w \rVert 
= \lim_{j\to\infty} \lVert Au - Au_j \rVert 
\le \alpha \lim_{j\to \infty} \lVert u - u_j\rVert
= 0.
\]
That is, $Au = w$. 
\end{proof}

A close read of this proof snippet
suggests that it doesn't seem to use the Hilbert space structure very much.
And indeed, it is possible to extract the work as a lemma in general metric spaces.
The appropriate abstraction is the following property of a function $f:X\to Y$ between metric spaces:
that there exists a constant $\beta>0$ such that for all $x_1$ and $x_2$,
$\beta d_X(x_1,x_2)\le d_Y(f(x_1),f(x_2))$.

As it turns out,
the same argument appears in the proof of the Contraction Mapping Theorem,
in a different special case (the case $Y = X$).
When Yury Kudryashov formalised the Contraction Mapping Theorem for Mathlib in 2020,
he recognised the appropriate context for the argument,\footnote{
  \url{https://github.com/leanprover-community/mathlib/pull/1859\#discussion\_r365490281}}
and wrote a self-contained theory development in Mathlib for such functions,\footnote{
  Mathlib \cite{Mathlib}, \href{https://github.com/leanprover-community/mathlib4/blob/21b3a48ba0f69a4ae9eb44dbacbc4de43a6290bf/Mathlib/Topology/MetricSpace/Antilipschitz.lean}{\texttt{Topology/MetricSpace/Antilipschitz}}}
for which he introduced the name \emph{antilipschitz maps}.
In fact, I would not be surprised to learn that this fairly short (1-2 pages of text) and easy theory
has been rediscovered and redeveloped many times, under many names.

With that abstraction and theory development available,
the snippet of the Lax--Milgram theorem we are discussing reduces simply to the following:

\begin{proof}
  $A$ is uniformly continuous and by \cref{eq:A-antilipschitz} it is antilipschitz, so it is injective and has closed range.
\end{proof}

Daniel Roca Gonz\'alez contributed this efficient proof of the Lax--Milgram theorem to Mathlib in 2022.\footnote{
  Mathlib \cite{Mathlib}, \href{https://github.com/leanprover-community/mathlib4/blob/21b3a48ba0f69a4ae9eb44dbacbc4de43a6290bf/Mathlib/Analysis/InnerProductSpace/LaxMilgram.lean}{\texttt{Analysis/InnerProductSpace/LaxMilgram}}}
(The theorem had earlier been formalised in Coq \cite{BCFMM17}, following a somewhat different proof.)

In this example we see illustrated Bourbaki's original argument in favour of abstraction: deduplication.
Formal mathematics is done at scale:
 it 
is written from the axioms up,
so nontrivial proofs form parts of a vast corpus;
writing formal mathematics is much more like writing an encyclopaedia than like writing an article.
At this scale,
the ``two, or several theories" united by an abstraction
are very likely all to turn up,
and simple efficiency favours using the abstraction.\footnote{
To believe that habitual abstraction really will avoid the repetition of proofs at large scale,
you must be something of a Platonist:
you must believe (as I do!) that the ``natural context'' of an argument is sufficiently unambiguous
that others who need it will formulate it in the same way,
and thus be led to stumble across your version.
}

\subsection{Smooth vector bundles}

My last example (a bit more technical than the others in this article)
is taken from the theory of smooth vector bundles in Lean,
which is joint work of mine with Floris van Doorn in 2022--23.

In this example, the particular definition of smooth vector bundle we chose for our theory matters.

\begin{definition}
A \emph{smooth vector bundle} with \emph{fibre} $F$ over a smooth manifold $B$ consists of
\begin{itemize}
  \item a collection of topological vector spaces indexed by $B$;
  \item a topology on the \emph{total space}, i.e.\ on their disjoint union;
  \item a collection of \emph{trivialisations}, each identifying the fibre-union over some open set
  $U\subseteq B$ homeomorphically with $U\times F$, commuting via projections with the identity on
  $U$, and fibrewise an isomorphism of topological vector spaces;
  \item with the property that for two trivialisations in the collection the induced map
  $U\cap V \to \operatorname{End}(F)$ is smooth.
\end{itemize}
\end{definition}

I will discuss two approaches to the proof of the following statement.

\begin{proposition}
  The total space of a smooth vector bundle is a smooth manifold.
\end{proposition}

Note that the fact that this is a theorem to be proved, rather than part of the definition,
is a consequence of our choice of definition.

Here is how you might prove this theorem in prose.
Since I didn't find a presentation of the theory of vector bundles in the literature
which started with precisely our definition,
this proof is not taken directly from real life.

\begin{proof}
Let $H$ be the model space for the smooth manifold $B$.  Given a trivialisation
$\psi =(\psi_b,\psi_f): \pi^{-1}(U)\to U \times F$ and a chart
$\varphi : V \mathrel{\overset{\sim}{\to}} \varphi(V)\subseteq H$ for $B$, define a candidate chart
\[
  \Phi_{\psi,\varphi}:\pi^{-1}(U\cap V) \to \varphi(U \cap V) \times F,
\]
\[
  \Phi_{\psi,\varphi}(p) := (\varphi(\psi_b(p)), \psi_f(p)).
\]
We need to check that for any two trivialisations $\psi_1, \psi_2$ and any two charts $\varphi_1$,
$\varphi_2$ the transition function $\Phi_{\psi_2,\varphi_2} \circ \Phi_{\psi_1,\varphi_1}^{-1}$ is
smooth.  This works out since $\psi_2\circ\psi_1^{-1}$, $\varphi_1$ and $\varphi_2$ are all smooth.
\end{proof}

Our Lean formalisation
uses Kobayashi--Nomizu's abstraction of a \emph{structure groupoid} \cite{KN63}
for a way in which a space is modelled on another space,
which is used there as the approach to the definition of smooth manifolds.
S\'ebastien Gou\"ezel developed this theory in Mathlib in 2019.\footnote{
  Mathlib \cite{Mathlib}, \href{https://github.com/leanprover-community/mathlib4/blob/21b3a48ba0f69a4ae9eb44dbacbc4de43a6290bf/Mathlib/Geometry/Manifold/ChartedSpace.lean}{\texttt{Geometry/Manifold/ChartedSpace}}
}

The advantage of our chosen definition of smooth vector bundle
is that, following a suggestion of Gou\"ezel,\footnote{
  Mathlib \cite{Mathlib}, \href{https://github.com/leanprover-community/mathlib4/blob/21b3a48ba0f69a4ae9eb44dbacbc4de43a6290bf/Mathlib/Geometry/Manifold/ChartedSpace.lean\#L139}{\texttt{Geometry/Manifold/ChartedSpace}, line 139}
} 
it too can be expressed using this structure groupoid abstraction.
  Let $H$ be the model space for the smooth manifold $B$.  Let $E$ be a smooth vector bundle over
  $B$ with fibre $F$. We consider the sequence
  \[
    E \quad \xrightarrow{\text{modelled on}} \quad B \times F \quad \xrightarrow{\text{modelled on}} \quad H \times F:
  \]
\begin{enumerate}
  \item $E$ is modelled on $B \times F$ with the charts being the trivialisations,
      and our vector bundle definition amounts to the condition that the transition functions between these charts lie in the \emph{smooth fibrewise-linear
       groupoid};
  \item $B \times F$ is in turn is modelled on $H \times F$ with the charts being the usual
    product manifold charts, and with the transition functions between these charts lying in the
    usual smooth manifold structure groupoid.
\end{enumerate}
In this language, here is an outline of our formalisation\footnote{
  Mathlib \cite{Mathlib}, \href{https://github.com/leanprover-community/mathlib4/blob/21b3a48ba0f69a4ae9eb44dbacbc4de43a6290bf/Mathlib/Geometry/Manifold/VectorBundle/Basic.lean\#L486}{\texttt{Geometry/Manifold/VectorBundle/Basic}, line 486}}
of the theorem.

\begin{proof}
  ``Modellings'' can be composed, so the modellings of $E$ on $B \times F$ and of $B \times F$ on
$H \times F$ yield a modelling of $E$ on $H \times F$. Structure groupoid properties can also
be composed, so the transition functions between these induced charts lie in the smooth manifold
structure groupoid for $B \times F$.
\end{proof}

This composition theorem for structure groupoids was formulated by us for the project;\footnote{
  Mathlib \cite{Mathlib}, \href{https://github.com/leanprover-community/mathlib4/blob/21b3a48ba0f69a4ae9eb44dbacbc4de43a6290bf/Mathlib/Geometry/Manifold/LocalInvariantProperties.lean\#L698}{\texttt{Geometry/Manifold/LocalInvariantProperties}, line 698}
} 
to our knowledge it does not appear in the literature.

This hierarchy of undigested abstractions is certainly a more obscure approach to this material
than would be acceptable in a traditional prose presentation.
But it has a certain elegance,
and it brings organisational assistance:
some work can be done cleanly at high level,
and the more painful direct manipulation of partially-defined smooth functions
appears only when checking the various preconditions for the abstractions to apply.

This is very much a slogan of formalisation: that it
incentivises abstraction 
to cope with the demands of writing proofs in full detail \cite{CT24}.
Gonthier \cite{Gon08Notices} similarly notes that his formalisation of the Four-Colour Theorem produced several abstractions, ``new and rather elegant nuggets of mathematics,'' as a byproduct.

\section{Conclusion}

In this article I discuss only the question of,
given a fixed statement,
what constitutes a good proof (formal or informal) of that statement.
An orthogonal question
is how to best express the development of a whole mathematical theory.\footnote{
A crude analogy is to consider the statements of a mathematical theory as a digraph,
with edges denoting easy implications
(some implications are easy in both directions and their edges are bidirectional).
To design a mathematical theory development,
you must select a spanning tree for this digraph.
}
This is a big question and it has produced an interesting literature \cite{ACKMRS20,Gou22}.

I have argued that good computer-formalised writing differs from good prose writing in two aspects:
its incorporation of algorithms and of abstractions.
These two aspects have an interesting commonality:
in prose writing, both represent breaks in tone, or even in the very experience of reading---%
moments at which the reader is sent to a reference in order to read up on an unfamiliar abstraction,
or to her computer to study and run a piece of code.
But in formalised writing these are not breaks:
prerequisites, computation  
and main argument form an integrated whole.

Monta\~no \cite{Mon12} argues that we experience a proof as beautiful according to the narrative experience of reading it,
the ``quality of its storytelling.''
In formalised mathematical writing,
more kinds of thinking can be incorporated without causing breaks in the narrative flow.
Our storytelling will be all the richer in consequence.

\begin{credits}
  \subsubsection{\ackname} I am grateful to Isabelle Petersen for assistance in typesetting the notes for this article.
  I also thank the audience of the first version of this material
  (a talk at the 2023 workshop ``Machine Assisted Proof'' at the Institute for Pure and Applied Mathematics),
  whose stimulating comments helped to sharpen the argument,
  and Tom Hales, John Harrison and Kim Morrison for useful comments on a draft.
  
  \subsubsection{\discintname}
  The author has no competing interests to declare that are
  relevant to the content of this article.
\end{credits}

\bibliographystyle{splncs04}
\bibliography{refs}

\begin{thebibliography}{10}
\providecommand{\url}[1]{\texttt{#1}}
\providecommand{\urlprefix}{URL }
\providecommand{\doi}[1]{https://doi.org/#1}

\bibitem{ACKMRS20}
Affeldt, R., Cohen, C., Kerjean, M., Mahboubi, A., Rouhling, D., Sakaguchi, K.: Competing inheritance paths in dependent type theory: A case study in functional analysis. In: Automated Reasoning: 10th International Joint Conference, IJCAR 2020, Paris, France, July 1-4, 2020, Proceedings, Part II. pp. 3--20. Springer-Verlag, Berlin, Heidelberg (2020). \doi{10.1007/978-3-030-51054-1_1}

\bibitem{Arn98}
Arnold, V.I.: On teaching mathematics. Russ. Math. Surv.  \textbf{53}(1),  229--236 (1998). \doi{10.1070/RM1998v053n01ABEH000005}

\bibitem{BCFMM17}
Boldo, S., Cl\'{e}ment, F., Faissole, F., Martin, V., Mayero, M.: A {C}oq formal proof of the {L}ax-{M}ilgram theorem. In: Proceedings of the 6th ACM SIGPLAN Conference on Certified Programs and Proofs. pp. 79--89. CPP 2017, Association for Computing Machinery, New York, NY, USA (2017). \doi{10.1145/3018610.3018625}

\bibitem{Bon82}
Bonsall, F.F.: A down-to-earth view of mathematics. The American Mathematical Monthly  \textbf{89}(1),  8--15 (1982). \doi{10.2307/2320989}

\bibitem{Bou50}
Bourbaki, N.: The architecture of mathematics. The American Mathematical Monthly  \textbf{57}(4),  221--232 (1950). \doi{10.2307/2305937}

\bibitem{CW07}
Chaieb, A., Wenzel, M.: Context aware calculation and deduction. In: Proceedings of the 14th Symposium on Towards Mechanized Mathematical Assistants: 6th International Conference. p. 27–39. Calculemus '07 / MKM '07, Springer-Verlag, Berlin, Heidelberg (2007). \doi{10.1007/978-3-540-73086-6_3}

\bibitem{CT24}
Commelin, J., Topaz, A.: Abstraction boundaries and spec driven development in pure mathematics. Bull. Am. Math. Soc.  \textbf{61}(2),  241--255 (2024). \doi{10.1090/bull/1831}

\bibitem{Mathlib}
mathlib community, T.: The {L}ean mathematical library. In: Proceedings of the 9th ACM SIGPLAN International Conference on Certified Programs and Proofs. p. 367–381. CPP 2020, Association for Computing Machinery, New York, NY, USA (2020). \doi{10.1145/3372885.3373824}, \url{https://github.com/leanprover-community/mathlib4/tree/21b3a48}, {Note:} The URL given is that for a Mathlib version from May 6, 2024 on Lean version 4.8.0-rc1; references to Mathlib files in this article are to this version.

\bibitem{CBG08}
Conway, J.H., Burgiel, H., Goodman-Strauss, C.: The symmetries of things. A K Peters, Wellesley, MA (2008). \doi{10.1201/b21368}

\bibitem{Eck}
Eck, D.: Wallpaper symmetry sketchpad, \url{https://math.hws.edu/eck/js/symmetry/wallpaper.html}, accessed: March 25, 2024

\bibitem{Gon08Notices}
Gonthier, G.: Formal proof - the four color theorem. Notices Am. Math. Soc.  \textbf{55}(11),  1382--1393 (2008), \url{www.ams.org/notices/200811/tx081101382p.pdf}

\bibitem{Gon13}
Gonthier, G., Asperti, A., Avigad, J., Bertot, Y., Cohen, C., Garillot, F., Le~Roux, S., Mahboubi, A., O'Connor, R., Biha, S.O., Pasca, I., Rideau, L., Solovyev, A., Tassi, E., Th\'{e}ry, L.: A machine-checked proof of the odd order theorem. In: Proceedings of the 4th International Conference on Interactive Theorem Proving. pp. 163--179. ITP'13, Springer-Verlag, Berlin, Heidelberg (2013). \doi{10.1007/978-3-642-39634-2_14}

\bibitem{Gou22}
Gou\"{e}zel, S.: A formalization of the change of variables formula for integrals in mathlib. In: Intelligent Computer Mathematics: 15th International Conference, CICM 2022, Tbilisi, Georgia, September 19-23, 2022, Proceedings. pp. 3--18. Springer-Verlag, Berlin, Heidelberg (2022). \doi{10.1007/978-3-031-16681-5_1}

\bibitem{Hal17}
Hales, T., Adams, M., Bauer, G., Dang, T.D., Harrison, J., Hoang, L.T., Kaliszyk, C., Magron, V., McLaughlin, S., Nguyen, T.T., et~al.: A formal proof of the {K}epler conjecture. Forum of Mathematics, Pi  \textbf{5}, ~e2 (2017). \doi{10.1017/fmp.2017.1}

\bibitem{Hal70}
Halmos, P.R.: How to write mathematics. Enseign. Math. (2)  \textbf{16},  123--152 (1970)

\bibitem{Har40}
Hardy, G.H.: A mathematician's apology. Cambridge {University} {Press} (1940)

\bibitem{Har96}
Harrison, J.: Formalized mathematics. Technical Report~36, Turku Centre for Computer Science (TUCS), Lemmink{\"a}isenkatu 14 A, FIN-20520 Turku, Finland (1996), \url{http://www.cl.cam.ac.uk/~jrh13/papers/form-math3.html}

\bibitem{Har07}
Harrison, J.: Automating elementary number-theoretic proofs using {G}r{\"o}bner bases. In: Pfenning, F. (ed.) Proceedings of the 21st International Conference on Automated Deduction, CADE 21. Lecture Notes in Computer Science, vol.~4603, pp. 51--66. Springer-Verlag, Bremen, Germany (2007). \doi{10.1007/978-3-540-73595-3_5}

\bibitem{Har09}
Harrison, J.: Formalizing an analytic proof of the {P}rime {N}umber {T}heorem (dedicated to {M}ike {G}ordon on the occasion of his 60th birthday). Journal of Automated Reasoning  \textbf{43},  243--261 (2009). \doi{10.1007/s10817-009-9145-6}

\bibitem{How20}
Howard, P.: Lecture notes for {M}612: {Partial Differential Equations} (2020), \url{https://people.tamu.edu/~phoward/M612.html}

\bibitem{Imm2018}
Immler, F.: A verified {ODE} solver and the {L}orenz attractor. Journal of Automated Reasoning  \textbf{61}(1),  73--111 (2018). \doi{10.1007/s10817-017-9448-y}

\bibitem{IA14}
Inglis, M., Aberdein, A.: Beauty is not simplicity: An analysis of mathematicians' proof appraisals. Philosophia Mathematica  \textbf{23}(1),  87--109 (07 2014). \doi{10.1093/philmat/nku014}

\bibitem{KN63}
Kobayashi, S., Nomizu, K.: Foundations of differential geometry. {I}, Intersci. Tracts Pure Appl. Math., vol.~15. Interscience Publishers, New York, NY (1963)

\bibitem{KS67}
Kochen, S., Specker, E.P.: The problem of hidden variables in quantum mechanics. Journal of Mathematics and Mechanics  \textbf{17}(1),  59--87 (1967)

\bibitem{Mon12}
Monta{\~n}o, U.: Ugly mathematics: why do mathematicians dislike computer-assisted proofs? Math. Intell.  \textbf{34}(4),  21--28 (2012). \doi{10.1007/s00283-012-9325-9}

\bibitem{Lean4}
Moura, L.d., Ullrich, S.: The {L}ean 4 theorem prover and programming language. In: Automated Deduction – CADE 28: 28th International Conference on Automated Deduction, Virtual Event, July 12–15, 2021, Proceedings. p. 625–635. Springer-Verlag, Berlin, Heidelberg (2021). \doi{10.1007/978-3-030-79876-5_37}

\bibitem{Lean3}
Moura, L.M.d., Kong, S., Avigad, J., van Doorn, F., von Raumer, J.: The {L}ean theorem prover (system description). In: Felty, A.P., Middeldorp, A. (eds.) Automated Deduction - {CADE-25} - 25th International Conference on Automated Deduction, Berlin, Germany, August 1-7, 2015, Proceedings. Lecture Notes in Computer Science, vol.~9195, pp. 378--388. Springer (2015). \doi{10.1007/978-3-319-21401-6_26}

\bibitem{Per91}
Peres, A.: Two simple proofs of the {Kochen}-{Specker} theorem. J. Phys. A, Math. Gen.  \textbf{24}(4),  l175--l178 (1991). \doi{10.1088/0305-4470/24/4/003}

\bibitem{Pot08}
Pottier, L.: Connecting {G}r{\"{o}}bner bases programs with {C}oq to do proofs in algebra, geometry and arithmetics. In: Rudnicki, P., Sutcliffe, G., Konev, B., Schmidt, R.A., Schulz, S. (eds.) Proceedings of the {LPAR} 2008 Workshops, Knowledge Exchange: Automated Provers and Proof Assistants, and the 7th International Workshop on the Implementation of Logics, Doha, Qatar, November 22, 2008. {CEUR} Workshop Proceedings, vol.~418. CEUR-WS.org (2008), \url{https://ceur-ws.org/Vol-418/paper5.pdf}

\bibitem{Rot97}
Rota, G.C.: The phenomenology of mathematical beauty. Synthese  \textbf{111}(2),  171--182 (1997). \doi{10.1023/A:1004930722234}

\bibitem{Tao07}
Tao, T.: What is good mathematics? Bull. Am. Math. Soc., New Ser.  \textbf{44}(4),  623--634 (2007). \doi{10.1090/S0273-0979-07-01168-8}

\bibitem{Wel90}
Wells, D.: Are these the most beautiful? Math. Intell.  \textbf{12}(3),  37--41 (1990). \doi{10.1007/BF03024015}

\end{thebibliography}

\end{document}